\documentstyle[amsfonts,12pt]{article} 

\newtheorem{theorem}{Theorem}

\begin{document}
\title{How weak is weak extent?}
\author{M. V. Matveev\\
{\small Department of mathematics, University of California at Davis,}\\
{\small Davis, CA 95616, USA (address valid till June 30, 2000)}\\
{\small {misha$\underline{\mbox{ }}$matveev@hotmail.com}}}
\date{}
\maketitle

{\bf Abstract.}{\small We show that the extent of a Tychonoff
space of countable weak extent can be arbitrary big.
The extent of $X$ is
$e(X)={\rm sup}\{|F|:F\subset X$ is closed and discrete$\}$
while
$we(X)={\rm min}\{\tau:\mbox{ for every open cover }{\cal U}
\mbox{ of }X\mbox{ there is }A\subset X\mbox{ such that }|A|\leq\tau
\mbox{ and }St(A,{\cal U})=X\}$ is the weak extent of $X$
(also called the star-Lindel\"of number of $X$).
Also we show that the extent of a normal space with countable
weak extent is not greater than $\bf c$.
}



\bigskip
Recall that the extent of a topological space $X$ is the cardinal
$e(X)={\rm sup}\{|F|:F\subset X$ is closed and discrete$\}$.
The {\em weak extent} of $X$ is the cardinal
$we(X)={\rm min}\{\tau:$ for every open cover $\cal U$ 
of $X$ there is $A\subset X$ such that $|A|\leq\tau$ 
and $St(A,{\cal U})=X\}$ \cite{Hod}.
The reason for this name is that for any $X\in{\rm T}_1$,
$we(X)\leq e(X)$; indeed, supposing $we(X)>\kappa$, there is an open cover
such that for every $A\subset X$ with $|A|\leq\kappa$ one has
$St(A,{\cal U})\neq X$; then one can inductively choose points
$x_\alpha$, $\alpha<\kappa$, so that 
$x_\alpha\notin St(\{x_\beta:\beta<\alpha\},{\cal U})$ for each $\alpha$;
once the points have been choosen the set 
$\{x_\alpha:\alpha<\kappa\}$ is closed, discrete
and of cardinality $\kappa$, so $e(X)\geq\kappa$.
Note also that $we(X)\leq d(X)$ obviously holds for every $X$.
Some cardinal inequalities involving extent can be improved by
replacing extent by weak extent.
Thus for $X\in T_1$, $|K(X)|\leq we(X)^{psw(X)}$ \cite{Hod}.
A natural question was stated in \cite{Milena}, \cite{Matsur}: 
how big can the difference between
the extent and the weak extent of a ${\rm T}_i$ space be?
First, we give the answer for the Tychonoff case.

\begin{theorem}\label{odin}
For every cardinal $\tau$ there is a Tychonoff space $X$
such that $e(X)\geq\tau$ and $we(X)=\omega$.
\end{theorem}

Before the paper \cite{Hod} the cardinal function $we(X)$
was called the {\em star-Lindel\"of number}
\cite{MM}, \cite{Matthes}, \cite{SW}.
In particular, a space $X$ such that $we(X)=\omega$ is called
star-Lindel\"of or $^*$Lindel\"of, see e.g. \cite{Ikenaga},
\cite{Milena}, \cite{Versus} \cite{Products}.

Note that $e(X)\leq 2^{we(X)\chi(X)}$ for every regular space $X$
\cite{Milena}.

To prove Theorem~1, we use a set-theoretic fact in Theorem~2 below..
Let $S$ be a set and $\lambda$ a cardinal.
A {\em set mapping} of order $\lambda$ on $S$ is a mapping that assigns to
each $s\in S$ a subset $f(s)\subset S$ so that 
$|f(s)|<\lambda$ and $s\not\in f(s)$.
A subset $T\subset S$ is called {\em $f$-free} if $f(t)\cap T=\emptyset$
for every $t\in T$. 
Answering a question of Erd\"os, Fodor proved in 1952
(\cite{Fodor}, see also \cite{Williams}, Theorem~3.1.5)
a general theorem a partial case of which is the following

\begin{theorem}
{\rm (Fodor)}
Let $S$ be a set of cardinality $\tau$ and let $f$ be a set mapping on $S$
of order $\omega$.
Then there is a countable family $\cal H$ of $f$-free subsets 
of $S$ such that $\cup{\cal H}=S$.
\end{theorem}

{\bf Proof of Theorem 1:}
Let $\tau$ be an infinite cardinal.
For each $\alpha<\tau$, $z_\alpha$ denotes the point in $D^\tau$
with only the $\alpha$-th coordinate equal to $1$.
Put $Z=\{z_\alpha:\alpha<\tau\}$.
Then $Z$ is a discrete subspace of $D^\tau$.
Further, let $\kappa$ be a cardinal such that 
${\rm cf}(\kappa)>\tau$. Put
$$
X=\left(D^\tau\times(k+1)\right)\setminus
\left(\left(D^\tau\setminus Z\right)\times\{\kappa\}\right).
$$
Also we denote $X_0=D^\tau\times\kappa$
and $X_1=Z\times\{\kappa\}=
\{(z_\alpha,\kappa):\alpha<\tau\}$.
Then $X=X_0\cup X_1$.

It is clear that $X_1$ is closed in $X$ and discrete, so $e(X)\geq\tau$.

It remains to prove that $we(X)=\omega$.
First, note that $X_0$ is countably compact, hence star-Lindel\"of.
So it remains to prove that $X_1$ is relatively star-Lindel\"of in $X$,
i.e. for every open cover $\cal U$ of $X$ there is a countable
$A\subset X$ such that $St(A,{\cal U}\}\supset X_1$.

Let $\cal U$ be an open cover of $X$.
For every $\alpha<\tau$ choose an $U_\alpha\in\cal U$ so that
$(z_\alpha,\kappa)\in U_\alpha$.
Further, for every $\alpha<\tau$ choose $\xi_\alpha<\kappa$
and $B_\alpha$, an element of the standard base of $D^\tau$,
so that $(z_\alpha,\kappa)\in(B_\alpha\times(\xi_\alpha,\kappa])\cap X
\subset U_\alpha$.
It remains to check that

\medskip
($+$) there is a countable $C\subset D^\tau$ such that
$B_\alpha\cap C\neq\emptyset$ for every $\alpha<\tau$.

\medskip\noindent
Indeed, since ${\rm cf}(\kappa)>\tau$, there is a $\gamma<\kappa$
such that $\gamma>\xi_\alpha$ for all $\alpha<\tau$.
Put $A=C\times\{\gamma\}$.
Then $U_\alpha\cap A\neq\emptyset$ for all $\alpha<\tau$,
so $X_1\subset St(A,{\cal U})$.
 
Now we check ($+$). The set $B_\alpha$ has the form
$$
B_\alpha=\{x\in D^\tau:x(\alpha)=1\mbox{ and }x(\alpha')=0
\mbox{ for all }\alpha'\in A_\alpha\}
$$
where $A_\alpha$ is some finite subset of $\tau\setminus\{\alpha\}$.
Consider the set mapping $f$ that assigns $A_\alpha$ to $\alpha$
for each $\alpha<\tau$.
By Fodor's theorem, there is a countable, $f$-free family
${\cal H}=\{H_n:n\in\omega\}$ of subsets of $\tau$
such that $\cup{\cal H}=\tau$.
For each $n\in \omega$, denote by $c_n$ the indicator function
of $H_n$, i.e. $c_n(\alpha)=1$ iff $\alpha\in H_n$.
Since $H_n$ is $f$-free, $B_\alpha\ni c_n$ for all $\alpha\in H_n$.
Put $C=\{c_n:n\in\omega\}$.
Then $B_\alpha\cap C\neq\emptyset$ for every $\alpha<\tau$,
i.e. ($+$) holds.

Pseudocompactness of $X$ follows from the fact that $X$ contains a dense countably
compact subspace $X_0$.
$\Box$

\medskip
Now we are going to show that in the normal case the extent of a space
of countable weak extent is not greater than $\bf c$.
In fact, we will prove a slightly more general statement.
Recall that a family of sets is linked if every two elements
have nonempty intersection.
The linked-Lindel\"of number of $X$ is the cardinal
$ll(X)=\min\{\tau:$ every open cover of $X$ has a subcover 
representable as the union of at most $\tau$ many linked
subfamilies$\}$ \cite{Products}.
A space $X$ with $ll(X)=\omega$ is called linked-Lindel\"of
\cite{Products}.
It is easy to see that $ll(X)\leq we(X)$ for every $X$.

\begin{theorem}\label{aboutNormal}
For every normal space $X$, $e(X)\leq 2^{ll(X)}$.
\end{theorem} 

{\bf Proof:}
Let $\tau$ be an infinite cardinal,
$K$ a closed discrete subspace of a normal space $X$
and $|K|=k>2^\tau$.
We have to show that $ll(X)>\tau$.
It is easy to construct a family $\cal A$ of subsets of $K$
such that $|{\cal A}|=k$ and for every nonempty finite
subfamily of $\cal A$, say $A_1,\dots,A_n,A_{n+1},\dots,A_{n+m}$,
$$
(*)\quad |A_1\cap\dots\cap A_n\cap (K\setminus A_{n+1})\cap\dots
\cap(K\setminus A_{n+m})|=k.
$$
For every $A\in\cal A$ pick a continuous function $f_A:X\to I$
such that $f(A)=\{1\}$ and $f(K\setminus A)=\{0\}$.
Denote ${\cal F}=\{f_A:A\in{\cal A}\}$ and
$F=\Delta{\cal F}:X\to I^{\cal F}$.
Then $|{\cal F}|=k$.
Note that $F(K)\subset D^{\cal F}$.
It follows from $(*)$ that $F(K)$ is dense in $D^{\cal F}$,
moreover, every open set in $D^{\cal F}$ containes $k$
elements of $F(K)$.
There is therefore a bijection $\varphi:K\to\cal B$,
where $\cal B$ is the standard base of $D^{\cal F}$,
such that $\varphi(z)\ni F(z)$ for every $z\in K$.
Every element $B\in\cal B$ has the form 
$$
B=B_{f_1\dots f_n}^{i_1\dots i_n}=\{x\in D^{\cal F}:
x(f_1)=i_1,\dots,x(f_n)=i_n\}
$$
where $n\in\Bbb N$, $f_1,\dots f_n\in\cal F$ and
$i_1,\dots, i_n\in D$.
Denote
$$
U(B)=\left\{x\in I^{\cal F}:\forall j\in\{1,\dots, n\}
\left(\begin{array}{c}
x(f_j)>\frac{1}{2}\mbox{ if }i_j=1\\
x(f_j)<\frac{1}{2}\mbox{ if }i_j=0
\end{array}\right)\right\}.
$$
Further, for every $z\in K$ put 
$\tilde{\varphi}(z)=U(\varphi^{-1}(z))$.
Then $\tilde{\varphi}(z)$ is a neighbourhood of $F(z)$ in $I^{\cal F}$.
Note that 
$$
(**)\quad \tilde{\varphi}(z)\cap\tilde{\varphi}(z')\neq\emptyset
\quad\mbox{  iff  }\quad
\varphi(z)\cap\varphi(z')\neq\emptyset
$$
Let $\cal G$ denote the family of all continuous functions form 
$X$ to $I$, $G=\Delta{\cal G}:X\to I^{\cal G}$,
$\pi:I^{\cal G}\to I^{\cal F}$ is the natural projection.
For each $z\in K$ denote
$\tilde{\tilde{\varphi}}(z)=\pi^{-1}({\varphi}(z))$.
Then $\tilde{\tilde{\varphi}}(z)$ is a neighbourhood of $G(z)$
in $I^{\cal G}$ and 
$$
(***)\quad \tilde{\tilde{\varphi}}(z)\cap\tilde{\tilde{\varphi}}(z)
\neq\emptyset\quad\mbox{  iff  }\quad
\tilde{\varphi}(z)\cap\tilde{\varphi}(z')\neq\emptyset.
$$
Last, for every $z\in K$ put 
$\tilde{\tilde{\tilde{\varphi}}}(z)=
(\tilde{\tilde{\varphi}}(z)\setminus G(K\setminus\{z\}))\cap G(X)$.
Then $\tilde{\tilde{\tilde{\varphi}}}(z)$ is a neighbourhood of 
$G(z)$ in $G(X)$ and
$$
(*{\rm v})\quad
\tilde{\tilde{\tilde{\varphi}}}(z)\cap\tilde{\tilde{\tilde{\varphi}}}(z')
\neq\emptyset\quad\mbox{  iff  }\quad
\tilde{\tilde{\varphi}}(z)\cap\tilde{\tilde{\varphi}}(z)\neq\emptyset.
$$
Put ${\cal U}_0=\{\tilde{\tilde{\varphi}}(z):z\in K\}$.
Since $G$ is a homeomorphic embedding,
$G(K)$ is closed in $G(X)$, so $O=G(X)\setminus G(K)$ is open
and hence ${\cal U}={\cal U}_0\cup\{O\}$
is an open cover of $G(X)$.

Since $w(D^{\cal F})>2^\tau$, $\cal B$, a base of $D^{\cal F}$,
is not representable as the union of at most $\tau$ many linked
subfamilies (see e.g. \cite{LM}).
By $(**)$, $(***)$ and $(*{\rm v})$ the same can be said
about the family ${\cal U}_0$.
Note that for every $z\in K$,
$\tilde{\tilde{\tilde{\varphi}}}(z)$ is the only element of $\cal U$
that containes $z$.
So $\cal U$ does not have a subcover representable as the union 
of at most $\tau$ many linked subfamilies
and thus $ll(X)=ll(G(X))>\tau$.$\Box$

\medskip
It is not clear whether the inequality in the previous theorem can be made
strict, even with star-Lindel\"of number instead of linked-Lindel\"of.

\bigskip
{\bf Asknowlegement.}
The author expresses his gratitude to Angelo Bella and to
Marion Scheepers for usefull discussions, in particular,
Angelo Bella has drawn authors attention to the paper \cite{Hod} 
and Marion Scheepers
has drawn author's attention to the book \cite{Williams}.

The paper was written while the author was visiting the University
of California, Davis. The author expresses his gratitude to
colleagues from UC Davis for their kind hospitality.



\begin{thebibliography}{9}

\bibitem{Milena}
M. Bonanzinga,
{\em Star-Lindel\"of and absolutely star-Lindel\"of spaces},
Q$\&$A in General Topology, {\bf 16} (1998) 79-104.

\bibitem{Versus}
M. Bonanzinga and M. V. Matveev,
{\em Star-Lindel\"ofness versus centered-Lindel\"ofness},
to appear in Comment. Math. Univ. Carol.

\bibitem{Products}
M. Bonanzinga and M. V. Matveev,
{\em Products of star-Lindel\"of and related spaces},
to appear in Houston J. of Math.

\bibitem{Fodor}
G. Fodor,
{\em Proof of a conjecture of P.~Erd\"os},
Acta Sci. Math. Szeged {\bf 14} (1952) 219-227.

\bibitem{Hod}
R. E. Hodel,
{\em Combinatorial set theory and cardinal function inequalities},
Proc. Amer. Math. Soc. {\bf 111} (1991) 567-575.

\bibitem{Ikenaga}
S. Ikenaga,
{\em A class which contains Lindel\"of spaces, 
separable spaces and countably compact spaces},
Memoires of Numazu College of Technology
{\bf 18} (1983) 105-108.

\bibitem{LM}
R. Levy and M.V. Matveev,
{\em Spaces with $\sigma$-$n$-linked topologies as special
subspaces of separable spaces},
Comment. Math. Univ. Carol. {\bf 40} (1999) 561-570.

\bibitem{Matthes}
M. V. Matveev,
{\em Pseudocompact and Related Spaces},
thesis, Moscow State University, Moscow, 1984.

\bibitem{Matsur}
M. V. Matveev,
{\em A survey on star covering properties},
Topological Atlas, Preprint No 330,\\
http://www.unipissing.ca/topology/v/a/a/a/19.htm

\bibitem{MM}
Dai MuMing,
{\em A topological space cardinality inequality involving the 
$^*$Lindel\"of number},
Acta Math. Sinica, {\bf 26} (1983) 731-735.

\bibitem{SW}
S. H. Sun and Y. M. Wang,
{\em A strengthened topological cardinal inequality},
Bull. Austral. Math. Soc. {\bf 32} (1985) 375-378.

\bibitem{Williams}
N.H.  Williams
{\em Combinatorial Set Theory},
North-Holland 1977.

\end{thebibliography}
\end{document}